\documentclass[final]{amsart}

\input xypic
\xyoption{all}

\usepackage{tikz-cd}
\usepackage{amsmath, amssymb} 
\usepackage{enumitem} %change enumerate labels

\usepackage{ifdraft}
\usepackage{todonotes}
\usepackage[final,colorlinks,linkcolor=blue!50!black,citecolor=green!50!black]{hyperref}

\newtheorem{theorem}{Theorem}[section]
\newtheorem{proposition}[theorem]{Proposition}
\newtheorem{lemma}[theorem]{Lemma}

\theoremstyle{definition}
\newtheorem{definition}[theorem]{Definition}
\newtheorem{remark}[theorem]{Remark}
\newtheorem{example}[theorem]{Example}

\newtheorem{problem}[theorem]{Problem}

\newcommand{\cohlgy}[1]{\ensuremath{H^{*}(#1)}} 

\newcounter{bean}

\newcommand{\namedright}[3]{\ensuremath{#1\stackrel{#2}
 {\longrightarrow}#3}}
\newcommand{\nameddright}[5]{\ensuremath{#1\stackrel{#2}
 {\longrightarrow}#3\stackrel{#4}{\longrightarrow}#5}}

\newcommand{\larrow}{\relbar\!\!\relbar\!\!\rightarrow} 
\newcommand{\llarrow}{\relbar\!\!\relbar\!\!\larrow}

\newcommand{\Q}{\mathbb{Q}}
\newcommand{\tensor}{\otimes}

\DeclareMathOperator{\Alt}{Alt}
\newcommand{\GL}{\mathrm{GL}}
\newcommand{\SL}{\mathrm{SL}}
\newcommand{\iso}{\cong}
\newcommand{\cf}{\mathfrak{c}}
\newcommand{\du}{\underline{d}}
\newcommand{\CAT}{\mathsf{CAT}}
\newcommand{\TOP}{\mathsf{TOP}}
\newcommand{\cell}{\mathit{cell}}

\newcommand{\Z}{\mathbb{Z}}
\newcommand{\N}{\mathbb{N}}

\begin{document} 

\title[Realizing orders]{Realizing orders in rational sphere product algebras with three generators}
\author{Tseleung So}
\address{Tseleung So, Institute of Mathematical Science, Pusan National University, Busan 46241, Republic of Korea}
\email{larry.so.tl@gmail.com}
\author{Donald Stanley}
\address{Donald Stanley, Department of Mathematics and Statistics, Univeristy of Regina, Regina, SK S4S 0A2, Canada}
\email{Donald.Stanley@uregina.ca}
\author{Stephen Theriault}
\address{Stephen Theriault, School of Mathematical Sciences, University of Southampton, Southampton SO17 1 BJ, United Kingdom}
\email{S.D.Theriault@soton.ac.uk}
\author{Ben Williams}
\address{Ben Williams, Department of Mathematics, The Univeristy of British Columbia, Vancouver, BC V6T 1Z2, Canada}
\thanks{Ben Williams acknowledges the support of the Natural Sciences and Engineering Research Council of Canada (NSERC), RGPIN-2021-02603. Donald Stanley acknowledges the support of the Natural Sciences and Engineering Research Council of Canada (NSERC), RGPIN-05466-2020. Donald Stanley \& Ben Williams acknowledge that this research was supported in part by the Pacific Institute for the Mathematical Sciences (PIMS) through CRG41.
Tseleung So aknowledges the support of the National Research Foundation of Korea (NRF) grant funded by the Korea government (MSIT) (RS-2025-0055914), and the support of NSERC Discovery Grant and NSERC RGPIN-2020-06428.
}
\email{tbjw@math.ubc.ca}

\subjclass[2020]{Primary 55N10 Secondary 13F55, 55P99, 55T20.} 
\keywords{graded order, cohomology realization problem, polyhedral product}

\begin{abstract}
The realization problem asks which algebras can be realized as the cohomology 
of spaces. We study this problem in the context of the orders in a graded 
rational exterior algebra on three generators. An order is a subring whose 
underlying additive group is a lattice. We give conditions for when such 
an order is realizable, and in particular show that in the simply-connected 
case any order is realizable if the generators of the exterior algebra are of odd degree. 
\end{abstract}

\maketitle 

\section{Introduction} 

The realization problem is a fundamental problem in algebraic topology. It asks: \emph{given an $S$-algebra $A$, can it be realized as the cohomology ring $A\cong H^*(X;S)$ of some space $X$?}
While it has been completely solved in certain cases, such as when $A$ is a simply connected rational algebra~\cite{quillen} or a integral polynomial ring~\cite{AG}, the general problem remains highly challenging and open.

Consider the realization problem for algebras having the form $A=F_{S}[x_1,\ldots,x_m]/I$, where $F_{S}[x_1,\ldots,x_m]$ is the free graded-commutative $S$-algebra on $m$ generators and $I$ is a homogeneous ideal.
The case $S=\Z$ has been studied extensively.
For example, when all generators have degree two and $I$ is generated by a finite number of square-free monic monomials, the resulting algebra is a \emph{Stanley-Reisner ring} which can be realized as the cohomology ring of a Davis-Janusziekwicz space~\cite{DJ}.
If the square-free condition on the generating monomials of $I$ is removed, then $A$ is a \emph{monomial ideal ring} and its realizability was proved independently by Bahri-Bendersky-Cohen-Gitler and Trevisan~\cite{BBCG2,trevisan}.
Further, if the generators $x_i$ are allowed to have arbitrary degrees and $I$ is generated by monic monomials, the first two authors of this paper showed that $A$ is isomorphic to the cohomology ring of some space modulo torsion~\cite{SS}.
Without modding out by torsion, $A$ is not realizable in general; 
explicit examples were given by Takeda~\cite{takeda}. 
%Suppose that all $x_i$ have even degrees, and that $x_ix_j=0$ whenever $|x_i|=|x_j|=2^r$ for some $r>1$. Takeda determined which Stanley-Reisner rings under these conditions are realizable~\cite{takeda}.
In a different direction, Takeda and the second author related the Stanley-Reisner case to a graph colouring problem~\cite{ST}.

When $S=\Q$ and $I=(x^2_1,\ldots,x^2_m)$, the algebra $F_{\mathbb{Q}}[x_{1},\ldots,x_{m}]/(x^2_1,\ldots,x^2_m)$ is the rational cohomology ring of a product of spheres. A subring $A$ of this algebra is an \emph{order} if its underlying additive group is a graded lattice. We ask whether such an order is realizable. In~\cite{SST} a family of these orders satisfying a certain combinatorial condition was shown to be realizable, which captured all orders if $m=2$ but missed many cases even when $m=3$. In this paper, we show that all orders are realizable when $m=3$ provided that each of $x_{1}$, $x_{2}$ and $x_{3}$ has degree~$\geq 2$ and no two of them are in the same even degree  
(see Theorem~\ref{th:mainCombined}). In particular, all simply-connected orders are realizable if the degrees of $x_{1}$, $x_{2}$ and $x_{3}$ are odd. 
\medskip 

To state our results precisely, some conventions and definitions are needed. 
%Let $\N$ denote the set of positive integers.
All rings and algebras in the paper are assumed to be unital and associative. The term ``$\Z$-algebra'' may be replaced by ``ring''. Rings and algebras are graded by nonnegative integers. If $A$ is a graded object and $m$ is an integer, we will write $A_m$ for the summand of degree $m$. If $a \in A_m$, we will also write $|a|=m$.

Given a commutative ground ring $S$, placed in degree $0$, the notation $F_S[x_1, \dots, x_r]$ denotes the free graded-commutative $S$-algebra generated by homogeneous elements $x_1, \dots, x_r$. The degrees of the elements $x_i$ are omitted from the notation, but will be indicated whenever necessary.

\subsection{Sphere-product algebras}\label{sec:sphere-prod-algebr}

Suppose $m$ is a positive integer. An \emph{$m$-generated sphere-product $S$-algebra} is an $S$-algebra $R$ that is isomorphic to the cohomology of a product of $m$ positive-dimensional spheres with coefficients in $S$. Algebraically, $R$ is a sphere-product $S$-algebra if it is isomorphic to 
\[ \frac{F_S[x_1, \dots, x_m]}{(x_1^2, \dots, x_m^2)}, \qquad |x_i| > 0 \quad \forall i \in \{1, \dots, m\}.\]
Observe that the list of degrees of the generators $d_i := |x_i|$, up to reordering, is an invariant of the isomorphism class of $R$ as a graded algebra.

In the case where all the $x_i$ are of odd degree, sphere-product algebras are exterior algebras $\Lambda_S(x_1, \dots, x_m)$. By contrast, if any $x_i$ is of even degree and $S$ is not of characteristic $2$, the algebra cannot be generated by anticommuting elements.

Suppose $R$ is a finite-dimensional graded $\Q$-algebra. A ring $A$ is an \emph{order} in $R$ if it has the following two properties:
\begin{itemize}
\item $A$ is a subring of $R$,
\item the underlying additive group of $A$ is generated by a homogeneous basis of the $\Q$-vector space $R$, that is, $A$ is a graded lattice in $R$.
\end{itemize}
If $A$ is a graded ring whose underlying abelian group is finitely generated and free, then $A$ is canonically embedded as an order in $\Q \tensor_\Z A$. In particular, a ring $A$ is an order in a sphere-product $\Q$-algebra if and only if $A$ is torsion-free and $\Q \tensor_\Z A$ is a sphere-product $\Q$-algebra.

The motivating problem of this paper is the following.
\begin{problem}\label{pr:mainProblem}
  Solve the realization problem for orders in sphere-product $\Q$-algebras.
\end{problem}
This problem includes, as a notable special case, solving the realization problem for torsion-free rings $A$ for which $\Q \tensor_\Z A$ is an exterior algebra. If the number of generators is $1$ then $\mathbb{Q}\otimes_{\mathbb{Z}} A$ is the rational cohomology of a sphere and if the number of generators is $2$ then Problem \ref{pr:mainProblem} is solved by~\cite[Remark 8.7]{SST}. In this paper we focus on the case 
$m=3$, which is significantly harder than the $m=2$ case.

\begin{theorem}[Main theorem] \label{th:mainCombined}
  Suppose $\du=(d_1,d_2,d_3)$ is a list of positive integers greater than~$1$ in which any repeated term is odd. Let $R$ be the sphere-product $\Q$-algebra generated by $x_1, x_2, x_3$ in degrees $\du=(d_1,d_2,d_3)$. If $A$ is an order in $R$, then there exists a $CW$-complex $X_A$ such that $H^*(X_A;\Z)\cong A$.
\end{theorem}

The proof is an immediate consequence of Theorems \ref{thm_top_main} and \ref{th:ordersAllBasisG} below. In order to state these results, we introduce 
the notion of an ``$m$-generated weighted sphere-product ring,'' which is a kind of order in a sphere-product $\Q$-algebra.

Let $\mathcal{P}([m])$ be the power set on $[m]=\{1,\ldots,m\}$  
and let $\sigma=\{i_{1},\ldots,i_{k}\}$ be an ordered subset of~$[m]$. 
A \emph{coefficient sequence} is a function $\mathfrak{c}: \mathcal{P}([m]) \to \N$, whose value at $\sigma \subseteq [m]$ is written $\mathfrak c_\sigma$, that further satisfies: 
\begin{itemize}
\item $\mathfrak c_\sigma = 1$ if $\sigma$ has no more than $1$ element;
\item $\mathfrak c_\sigma \mathfrak c_\tau \mid \mathfrak c_{\sigma \cup \tau}$ for all $\sigma, \tau \subset [m]$ satisfying $\sigma \cap \tau = \emptyset$.
\end{itemize} 
Let $\mathfrak c$ be a coefficient sequence and let $\du = (d_1, \dots d_m)$ be a list of positive integers. Consider the sphere-product $\Q$-algebra
\[R=  \frac{F_\Q[a_1, \dots, a_m]}{(a_1^2, \dots, a_m^2)}, \quad |a_i| = d_i.\]

\begin{definition}\label{def:mgwspr}
  The \emph{($m$-generated) weighted sphere-product ring} $A(\mathfrak c, \du)$ is the order in $R$ that, as a lattice, is generated by the elements
\[ a_\sigma := \frac{1}{c_\sigma} \prod_{i \in \sigma} a_i \quad \text{(the product is taken in ascending order of $i$)}, \]
where $\sigma$ ranges over all elements in $\mathcal{P}([m])$.
\end{definition}

Our conventions imply that $a_\emptyset = 1$ and $a_{\{i\}} = a_i$. There are product relations:
\begin{equation}\label{dfn_A(c) prod}
a_\sigma a_\tau =
\begin{cases}
 \pm \dfrac{\cf_{\sigma \cup \tau}}{\cf_\sigma \cf_\tau} a_{\sigma \cup \tau} \qquad &\text{if $\sigma \cap \tau = \emptyset$;}\\
  0 \qquad &\text{otherwise.}
\end{cases}
\end{equation}
The ambiguous sign $\pm$ that appears above depends on $\sigma$, $\tau$ and on the parity of the degrees $d_i$. A precise formula for this sign may be given: one counts the number of pairs $i < j$ where $j \in \sigma$ and $i \in \tau$ for which $d_i, d_j$ are both odd. If this count is odd, then the sign is negative, otherwise it is positive.

We will say that $A$ is an ``$m$-generated weighted sphere-product ring'' if it is isomorphic to $A(\cf, \du)$ for some choices of $\cf$ and $\du$. If we take such an $A$ and apply $\Q \tensor_\Z A$, then the result is a sphere-product $\Q$-algebra. We conclude that weighted sphere-product rings are orders in sphere-product $\Q$-algebras.

Theorem~\ref{th:mainCombined} follows from two theorems concerning weighted sphere-product rings in the case of $m=3$ generators. The first says that all such rings are realizable provided we avoid generators in degree $1$.

\begin{theorem}[Topological existence theorem]\label{thm_top_main}
Let $\underline{d}=(d_1,d_2,d_3)$ be a list of positive integers greater than $1$, and let $\mathfrak{c} : \mathcal{P}([3]) \to \N$ be a coefficient sequence. 
There exists a $CW$-complex $X(\cf,\du)$ such that $H^*(X(\cf,\du);\Z)\cong A(\cf,\du)$.
\end{theorem} 

The second theorem says that weighted sphere-product rings account for a great many of the orders in $3$-generated sphere-product $\Q$-algebras, and in particular, for all orders in $3$-generated exterior algebras.

\begin{theorem}[Algebraic classification theorem] \label{th:ordersAllBasisG}
  Suppose $\du=(d_1,d_2,d_3)$ is a list of positive integers, where any repeated term is odd. Let $R$ be the sphere-product $\Q$-algebra generated by $x_1, x_2, x_3$ in degrees $\du=(d_1,d_2,d_3)$. If $A$ is an order in $R$, then there exists a coefficient sequence $\cf$ for which $A \iso A(\cf, \du)$.
\end{theorem} 

The weighted polyhedral product, which was introduced in~\cite{SST}, is the main tool used to prove the topological realization result, Theorem~\ref{thm_top_main}. The algebraic result, Theorem~\ref{th:ordersAllBasisG}, is established by working with presentations for sphere-product $\mathbb{Q}$-algebras.

Two cases remain to obtain a complete solution to the $3$-generator case of Problem \ref{pr:mainProblem}: algebras in which there are generators in degree~$1$ and algebras in which there is more than one generator in a given even degree. The degree-$1$ case is unsolved because weighted polyhedral products are most developed in the simply-connected case (see Theorem~\ref{realization}, for example). The case of more than one generator in a given even remains open because there are examples (see Example~\ref{ex:bad3}) of orders in weighted sphere-product $\mathbb{Q}$-algebras that are not weighted sphere-product rings.

The authors thank the anonymous referee for carefully reading the paper and for providing helpful comments.

\section{Weighted Polyhedral products}
Throughout this paper we assume that all spaces are pointed $CW$-complexes and all maps are pointed cellular. For $\{i_1,\ldots,i_k\}\subseteq[m]$, let $d_{\{i_1,\ldots,i_k\}}=d_{i_1}+\cdots+d_{i_k}$. In particular, $d_{[m]}=d_{1}+\cdots + d_{m}$.

Recall that a relative $CW$-pair $(X,A)$ is a pair of pointed $CW$-complexes such that $A$ is a subcomplex of $X$ and the inclusion $A\hookrightarrow X$ is a pointed map. Moreover let $(Y,B)$ be another relative $CW$-pair. A map of relative $CW$-pairs $f:(X,A)\to(Y,B)$ is a pointed map $f:X\to Y$ satisfying the commutative diagram
\[
\xymatrix{
A\ar[r]^-{f|_A}\ar@{^{(}->}[d]	&B\ar@{^{(}->}[d]\\
X\ar[r]^-{f}					&Y
}
\]

\begin{definition}\label{dfn_poly prod}
Let $(\underline{X},\underline{A})=\{(X_i,A_i)\}^m_{i=1}$ be a sequence of $m$ relative $CW$-pairs. For each $\sigma\subseteq[m]$, define $(\underline{X},\underline{A})^{\sigma}$ as the subspace
\[
(\underline{X},\underline{A})^{\sigma}=\left\{(x_1,\ldots,x_m)\in\prod^m_{i=1}X_i\mid x_i\in A_i\text{ whenever }i\notin\sigma\right\}
\]
of $\prod^m_{i=1}X_i$. Let $K$ be a simplicial complex on $[m]$. The \emph{polyhedral product} $(\underline{X},\underline{A})^{K}$ is defined by
\[
(\underline{X},\underline{A})^{K}=\bigcup_{\sigma\in K}(\underline{X},\underline{A})^{\sigma}\subseteq\prod^m_{i=1}X_i.
\]
Note that $(\underline{X},\underline{A})^{K}=\mbox{colim}_{\sigma\in K}(\underline{X},\underline{A})^{\sigma}$.
\end{definition}

\begin{definition}\label{dfn_poly prod functoriality}
Let $(\underline{X},\underline{A})=\{(X_i,A_i)\}^m_{i=1}$ and $(\underline{Y},\underline{B})=\{(Y_i,B_i)\}^m_{i=1}$ be sequences of $m$ relative $CW$-pairs and let $\{f_i:(X_i,A_i)\to(Y_i,B_i)\}^m_{i=1}$ be a sequence of maps of relative $CW$-pairs. The \emph{induced map} $f(K):(\underline{X},\underline{A})^{K}\longrightarrow(\underline{Y},\underline{B})^{K}$ is given by the restriction of $\prod^m_{i=1}f_i:\prod^m_{i=1}X_i\to\prod^m_{i=1}Y_i$ to $(\underline{X},\underline{A})^K$.
\end{definition}

We give an example that plays an important role in the paper. This requires some notation. 
Let $I=[0,1]$ be the unit interval and let $X$ be a space. The \emph{reduced cone} on $X$ is the 
quotient space 
\[CX=(I\times X)/\sim\] 
where $(0,x)\sim(0,x')$ and $(t,\ast)\sim (0,\ast)$, and the \emph{reduced suspension} on $X$ is the quotient space 
\[\Sigma X=(I\times X)/\sim\] 
where $(0,x)\sim(0,x')$, $(1,x)\sim(1,x')$ and $(t,\ast)\sim (0,\ast)$.

\begin{example}\label{ex_higher wh prod}
Given $m$ positive integers $d_1,\ldots,d_m$, let
\[
(\underline{CS},\underline{S})=\{(CS^{d_i-1},S^{d_i-1})\}^m_{i=1}\quad
\text{and}\quad
(\underline{\Sigma S},\ast)=\{(S^{d_i},\ast)\}^m_{i=1}.
\]
Then $(\underline{CS},\underline{S})^{\partial\Delta^{m-1}}\cong S^{d_{[m]}-1}$ and $(\underline{\Sigma S},\ast)^{\partial\Delta^{m-1}}$ is the fat wedge
\[
FW(S^{d_1},\ldots,S^{d_m})=\left\{(x_1,\ldots,x_m)\in\prod^m_{i=1}S^{d_i}\mid\text{at least one of $x_i$'s is }\ast\right\}.
\]
For $1\leq i\leq m$ let $q_i:(CS^{d_{i}-1},S^{d_{i}-1})\to(S^{d_{i}},\ast)$ be the quotient map collapsing the boundary $S^{d_i-1}$ of $CS^{d_i-1}$ to a point. Then the induced map $q(\partial\Delta^{m-1}):S^{d_{[m]-1}}\to FW(S^{d_1},\ldots,S^{d_m})$ is a higher Whitehead product, its mapping cone is $(\underline{\Sigma S},\ast)^{\Delta^{m-1}}\cong\prod^m_{i=1}S^{d_i}$, and there is a pushout:
\[
\xymatrix 
{S^{d_{123}-1}\cong (\underline{CS},\underline{S})^{\partial\Delta^2} 
  \ar[rr] \ar[d] &&  (\underline{CS},\underline{S})^{\Delta^2}  \ar[d] \\
  (\underline{S},*)^{\partial\Delta^2} \ar[rr]  && (\underline{S},*)^{\Delta^2} 
 }
 \]
\end{example}

In~\cite{SST}, weighted polyhedral products are defined as a generalization of polyhedral products. The definition requires a bit of setup.

\begin{definition}\label{dfn_power couple}
A \emph{$CW$-power couple} is a pointed relative $CW$-complex $(X,A)$ equipped with a collection of maps of pointed relative $CW$-complexes $\{\rho_a:(X,A)\to(X,A)\}_{a\in\N}$, called \emph{power maps}, such that $\rho_1$ is the identity map and $\rho_a\circ\rho_b=\rho_{ab}$ for any $a,b\in\N$. We say that $f:(X,A)\to(Y,B)$ is a \emph{map between $CW$-power couples} if it is a map of pointed relative $CW$-complexes $f:(X,A)\to(Y,B)$ making the diagrams
\begin{equation}\label{diagram_maps between power couples}
\xymatrix{
A\ar[r]^-{f|_A}\ar[d]_-{\rho_a|_A}	&B\ar[d]^-{\widetilde{\rho}_a|_B}\\
A\ar[r]^-{f|_A}						&B
}
\qquad
\xymatrix{
X\ar[r]^-{f}\ar[d]_-{\rho_a}	&Y\ar[d]^-{\widetilde{\rho}_a}\\
X\ar[r]^-{f}					&Y
}
\end{equation}
commute for each $a\in Z_+$, where $\rho_a$ and $\widetilde{\rho}_a$ are power maps for $(X,A)$ and $(Y,B)$.
\end{definition}

\begin{example}\label{ex_reduced suspension}
Identify $S^1$ with the unit circle $\{e^{ti}\mid t\in\mathbb{R}\}$ in $\mathbb{C}$. For $a\in\N$ define the power map on $(S^1,\ast)$ to be $e^{ti}\mapsto e^{ati}$. Let $(X,A)$ be a relative $CW$-pair. The relative $CW$-pair $(S^1\wedge X,S^1\wedge A)$ is a $CW$-power couple where the power map $\rho_a$ is given by
\[
\rho_a(e^{ti}\wedge x)=e^{ati}\wedge x.
\]
In particular $(S^d,\ast)$ is a $CW$--power couple for any $d\in\N$.
\end{example}

\begin{example}\label{ex_cone}
Suppose $(X,\ast)$ is a $CW$-power couple with power maps $\rho'_a$. Then $(CX,X)$ and $(\Sigma X, \ast)$ are $CW$-power couples where the power maps $\rho_a:(CX_{i},X_{i})\to(CX_{i},X_{i})$ and $\widetilde{\rho}_a:(\Sigma X_{i},\ast)\to(\Sigma X_{i},\ast)$ are given by
\[
\rho_a([t,x])=[t,\rho'_a(x)]\quad
\text{and}\quad
\widetilde{\rho}_a([[t,x]])=[[t,\rho'_a(x)]]
\]
respectively, with $[t,x]\in CX$ and $[[t,x]]\in\Sigma X$ being the quotient images of $(t,x)\in I\times X$.

If we take $(X,\ast)$ to be $(S^d,\ast)$ defined in Example~\ref{ex_reduced suspension}, then $(CS^d,S^d)$ and $(\Sigma S^d,\ast)$ are $CW$-power couples. 
\end{example}

\begin{definition}\label{dfn_power seq}
Fix a positive integer $m$. A \emph{power sequence} $c$ is a map
\[
c:\Delta^{m-1}\to\N^m,\qquad
\sigma\mapsto(c^{\sigma}_1,\ldots,c^{\sigma}_m)
\]
such that $c^{\sigma}_i=1$ for $i\notin\sigma$ and $c^{\tau}_i$ divides $c^{\sigma}_i$ for $\tau\subseteq\sigma$ and $1\leq i\leq m$. Further, let $(\underline{X},\underline{A})=\{(X_i,A_i)\}^m_{i=1}$ be a sequence of $CW$--power couples and let $K$ be a simplicial complex on $[m]$. If $\tau$ and $\sigma$ are subsets of $[m]$ such that $\tau\subseteq\sigma$, then ${c^{\sigma}_i/c^{\tau}_i}$ %$\displaystyle\frac{c^{\sigma}_i}{c^{\tau}_i}$
is a positive integer for each $i$. Denote by
\[
\underline{c}^{\sigma/\tau}\colon(\underline{X},\underline{A})^K\to (\underline{X},\underline{A})^K
\]
the restriction of the product of power maps $\prod^m_{i=1}{\rho_{c^{\sigma}_i/c^{\tau}_i}}:\prod^m_{i=1}X_i\to \prod^m_{i=1}X_i$.

From a power sequence $c$ we define a \emph{coefficient sequence} by  $\Phi(c)_\sigma=\prod_i c^{\sigma}_i$
\end{definition}

For a simplicial complex $K$, let $\CAT(K)$ be the \emph{face category} of $K$, where the 
objects are the faces (simplices) of $K$ and the maps are the inclusions of a sub-face into a face. Let 
$\TOP_{\ast}$ be the category of pointed topological spaces and pointed continuous maps between them.

\begin{definition}\label{dfn_weighted poly prod system}
Given a positive integer $m$, a sequence $(\underline{X},\underline{A})=\{(X_i,A_i)\}^m_{i=1}$ of $CW$--power couples and a power sequence $c$, a \emph{weighted polyhedral product system} $\{(\underline{X},\underline{A})^{\bullet,c},\eta(\bullet)\}$ consists of a functor $(\underline{X},\underline{A})^{\bullet,c}:\CAT(\Delta^{m-1})\to \TOP_*$ and a pointed map $\eta(\sigma):(\underline{X},\underline{A})^{\sigma}\to(\underline{X},\underline{A})^{\sigma,c}$ for each simplex $\sigma\in\Delta^{m-1}$ 
satisfying: 
\begin{enumerate}[label=(\roman*)]
\item	$(\underline{X},\underline{A})^{\emptyset,c}=(\underline{X},\underline{A})^{\emptyset}=\prod^m_{i=1}A_i$ and $\eta(\emptyset):(\underline{X},\underline{A})^{\emptyset}\to(\underline{X},\underline{A})^{\emptyset,c}$ is the identity map;
\item	for $1\leq i\leq m$, the space $(\underline{X},\underline{A})^{\{i\},c}$ and the maps $\eta(\{i\})$ and $i_{\emptyset}^{\{i\}}$ are defined by the pushout 
\[
\xymatrix{
(\underline{X},\underline{A})^{\emptyset}\ar@{^{(}->}[r]\ar[d]^-{\underline{c}^{\{i\}/\emptyset}}	&(\underline{X},\underline{A})^{\{i\}}\ar[d]^-{\eta(\{i\})}\\
(\underline{X},\underline{A})^{\emptyset,c}\ar[r]^-{\imath^{\{i\}}_{\emptyset}}								&(\underline{X},\underline{A})^{\{i\},c}
}
\]
where $\underline{c}^{\{i\}/\emptyset}$ is the product map $1\times\cdots\times c^{\{i\}}_i\times\cdots\times 1:\prod^m_{j=1}A_j\to\prod^m_{j=1}A_j$;
\item\label{dfn_weighted poly prod po}
for $\sigma\in\Delta^{m-1}$ with $|\sigma|>1$, the space $(\underline{X},\underline{A})^{\sigma,c}$ and the maps $\eta(\sigma)$ and $i_{\partial\sigma}^{\sigma}$ are defined by the pushout 
\[
\xymatrix{
(\underline{X},\underline{A})^{\partial\sigma}\ar[r]^-{\jmath^{\sigma}_{\partial\sigma}}\ar[d]^-{\eta(\partial\sigma)}	&(\underline{X},\underline{A})^{\sigma}\ar[d]^-{\eta(\sigma)}\\
(\underline{X},\underline{A})^{\partial\sigma,c}\ar[r]^-{\imath^{\sigma}_{\partial\sigma}}		&(\underline{X},\underline{A})^{\sigma,c}
}
\]
where, if $\tau\hookrightarrow\sigma$ and $\imath^{\sigma}_{\tau}:(\underline{X},\underline{A})^{\tau,c}\to(\underline{X},\underline{A})^{\sigma,c}$ is the map determined by $(\underline{X},\underline{A})^{\bullet,c}$, then: 
\begin{enumerate}[label=(\alph*)] 
\item $\jmath^{\sigma}_{\partial\sigma}$ is the inclusion; 
\item	$(\underline{X},\underline{A})^{\partial\sigma,c}= \mbox{colim}_{\tau\subsetneq\sigma}(\underline{X},\underline{A})^{\tau,c}$ for $\sigma\in\Delta^{m-1}$;
%\item	$\imath^{\sigma}_{\partial\sigma}:(\underline{X},\underline{A})^{\partial\sigma,c}\to(\underline{X},\underline{A})^{\sigma,c}$ is the colimit of $\imath_{\tau}^{\sigma}$ over $\tau\subsetneq\sigma$; 
\item	$\eta(\partial\sigma):(\underline{X},\underline{A})^{\partial\sigma}\to(\underline{X},\underline{A})^{\partial\sigma,c}$ is the colimit of composites
\[
(\underline{X},\underline{A})^{\tau}\overset{\underline{c}^{\sigma/\tau}}{\longrightarrow}(\underline{X},\underline{A})^{\tau}\overset{\eta(\tau)}{\longrightarrow}(\underline{X},\underline{A})^{\tau,c}\overset{\imath_{\tau}^{\partial\sigma}}{\longrightarrow}(\underline{X},\underline{A})^{\partial\sigma,c}
\]
over $\tau\subsetneq\sigma$, for $\underline{c}^{\sigma/\tau}$ as in Definition~\ref{dfn_power seq} and 
$\imath^{\partial\sigma}_{\tau}$ the inclusion into the colimit. 
\end{enumerate}
\end{enumerate}
\end{definition}

Examples will appear momentarily; we first 
state the realization result in~\cite[Theorem 8.4]{SST}. 
 
 \begin{theorem}\label{realization}
 Let $(\underline{\Sigma S},\ast)=\{(S^{d_i},\ast)\}^m_{i=1}$ with the power couple structure of Example \ref{ex_reduced suspension}
and $c$ a power sequence, Then \[A(\Phi(c), \underline d)\cong\cohlgy{(\underline{\Sigma S},\underline{\ast})^{\Delta^{m-1},c}}\] 
 \end{theorem} 

This relates to sphere product algebras. By Definition~\ref{dfn_poly prod}, 
$(\Sigma S,\ast)^{\Delta^{m-1}}=S^{d_{1}}\times\cdots\times S^{d_{m}}$, 
and by~\cite[Theorem 6.5]{SST} there is an isomorphism of $\mathbb{Z}$-modules 
$\cohlgy{(\underline{\Sigma S},\underline{\ast})^{\Delta^{m-1},c}}\cong 
   \cohlgy{(\underline{\Sigma S},\underline{\ast})^{\Delta^{m-1}}}$, 
although the ring structure may be different. 

\begin{remark}
Although $(\underline{X},\underline{A})^{\bullet,c}$ is a functor, $\eta(\bullet)=\{\eta(\sigma)\mid\sigma\subseteq[m]\}$ is not a natural transformation. In fact for any $\tau\subseteq\sigma\subseteq[m]$ there is a commutative diagram
\[
\xymatrix{
(\underline{X},\underline{A})^{\tau}\ar[r]^-{\jmath^{\sigma}_{\tau}}\ar[d]^-{\underline{c}^{\sigma/\tau}}	&(\underline{X},\underline{A})^{\sigma}\ar[dd]^-{\eta(\sigma)}\\
(\underline{X},\underline{A})^{\tau}\ar[d]^-{\eta(\tau)}	&\\
(\underline{X},\underline{A})^{\tau,c}\ar[r]^-{\imath^{\sigma}_{\tau}}	&(\underline{X},\underline{A})^{\sigma,c}
}
\]
\end{remark}

\begin{definition}
Let $\{(\underline{X},\underline{A})^{\bullet,c},\eta(\bullet)\}$ be a weighted polyhedral product system, and let $K$ be a simplicial complex on $[m]$. The \emph{weighted polyhedral product} $(\underline{X},\underline{A})^{K,c}$ is the colimit
\[
(\underline{X},\underline{A})^{K,c}=\mbox{colim}_{\sigma\in K}(\underline{X},\underline{A})^{\sigma,c}
\]
and its \emph{associated map} $\eta(K):(\underline{X},\underline{A})^{K}\to(\underline{X},\underline{A})^{K,c}$ is the colimit of composites
\[
(\underline{X},\underline{A})^{\sigma}\overset{\underline{c}^{\omega/\sigma}}{\longrightarrow}(\underline{X},\underline{A})^{\sigma}\overset{\eta(\sigma)}{\longrightarrow}(\underline{X},\underline{A})^{\sigma,c}\overset{\imath^K_{\sigma}}{\longrightarrow}(\underline{X},\underline{A})^{K,c}
\]
over $\sigma\in K$, where $\omega$ is the minimal simplex in $[m]$ containing $K$. 
\end{definition} 

\begin{example}\label{eg_ctrivial} 
Suppose the power sequence $c$ is trivial, that is $c(\sigma)=(1,\ldots,1)$ for any $\sigma\subseteq[m]$. Then $(\underline{X},\underline{A})^{\sigma,c}=(\underline{X},\underline{A})^{\sigma}$ and $\eta(\sigma)$ is the identity map.
\end{example}

\begin{example}\label{eg_Moore sp}
Take $m=3$. Let $d_1,d_2$ be positive integers greater than 1, and let $(\underline{X},\underline{A})=\{(X_i,A_i)\}^3_{i=1}$ be a sequence of $CW$--power couples where
\[
(X_i,A_i)=\begin{cases}
(CS^{d_i-1},S^{d_i-1})	&i=1,2\\
(\ast,\ast)				&i=3.
\end{cases}
\]
Here, the power maps of $(CS^{d_i-1},S^{d_i-1})$ are described in Example~\ref{ex_reduced suspension}. Suppose $c:P([3])\to\N^3$ satisfies $c_{\emptyset}=c(\{i\})=(1,1,1)$ for $1\leq i\leq 3$. By Definition~\ref{dfn_weighted poly prod system} we have
\[
\begin{array}{l}
(\underline{X},\underline{A})^{\{1,2\},c}\cong P^{d_{\{1,2\}}}(c_1^{\{1,2\}}c_2^{\{1,2\}})\\[5pt]
(\underline{X},\underline{A})^{\{2,3\},c}\cong(\underline{X},\underline{A})^{\{2\},c}=S^{d_1-1}\times CS^{d_2-1}\\[5pt]
(\underline{X},\underline{A})^{\{1,3\},c}\cong(\underline{X},\underline{A})^{\{1\},c}=CS^{d_1-1}\times S^{d_2-1}
\end{array}
\]
where $P^d(k)$ is the mapping cone of the degree $k$ map $S^{d-1}\to S^{d-1}$. Hence we have
\[
(\underline{X},\underline{A})^{\partial\Delta^2,c}=\mbox{colim}_{\{i_1,i_2\}\subset[3]}(\underline{X},\underline{A})^{\{i_1,i_2\},c}\cong {P^{d_{\{1,2\}}}(c^{\{1,2\}}_1c^{\{1,2\}}_2)}.
\]
\end{example}

As for polyhedral products in Definition~\ref{dfn_poly prod functoriality} we wish to construct a map $(\underline{X},\underline{A})^{\sigma,c}\to(\underline{Y},\underline{B})^{\sigma,c}$ from a collection of maps between $CW$--power couples $f_i:(X_i,A_i)\to(Y_i,B_i)$.

% \stanley{which proposition?}

\begin{proposition}[Proposition~4.1 in~\cite{SST}]\label{lemma_map between weight poly prod}
Let $c$ be a power sequence, let $(\underline{X},\underline{A})=\{(X_i,A_i)\}^m_{i=1}$ and $(\underline{Y},\underline{B})=\{(Y_i,B_i)\}^m_{i=1}$ be sequences of $CW$--power couples, and let
\[
\{(\underline{X},\underline{A})^{\bullet,c},\eta(\bullet)\}\quad
\text{and}\quad
\{(\underline{Y},\underline{B})^{\bullet,c},\widetilde{\eta}(\bullet)\}
\]
be the associated weighted polyhedral products systems. Suppose that $\{f_i:(X_i,A_i)\to(Y_i,B_i)\}^m_{i=1}$ is a sequence of maps between $CW$--power couples. Then there exists a unique collection of maps 
\[
f(\bullet,c)=
\{f(\sigma,c):(\underline{X},\underline{A})^{\sigma,c}\to(\underline{Y},\underline{B})^{\sigma,c}\mid\sigma\in\Delta^{m-1}\},
\] 
such that, for each $\tau\subseteq\sigma\in\Delta^{m-1}$, the following diagrams commute
\begin{enumerate}[label=(\alph*)]
\begin{minipage}{0.45\linewidth}
\item\label{diagram_unique weight poly prod1}
\hspace{0.5cm}
$\xymatrix{
(\underline{X},\underline{A})^{\tau,c}\ar[r]^-{f(\tau,c)}\ar[d]^-{\imath^{\sigma}_{\tau}}	&(\underline{Y},\underline{B})^{\tau,c}\ar[d]^-{\widetilde{\imath}^{\,\sigma}_{\tau}}\\
(\underline{X},\underline{A})^{\sigma,c}\ar[r]^-{f(\sigma,c)}		&(\underline{Y},\underline{B})^{\sigma,c}
}$
\end{minipage}
\hfill
\begin{minipage}{0.45\linewidth}
\item\label{diagram_unique weight poly prod2}
\hspace{0.5cm}
$\xymatrix{
(\underline{X},\underline{A})^{\sigma}\ar[r]^-{f(\sigma)}\ar[d]^-{\eta(\sigma)}	&(\underline{Y},\underline{B})^{\sigma}\ar[d]^-{\widetilde{\eta}(\sigma)}\\
(\underline{X},\underline{A})^{\sigma,c}\ar[r]^-{f(\sigma,c)}		&(\underline{Y},\underline{B})^{\sigma,c}
}$
\end{minipage}
\end{enumerate}
where $f(\sigma)$ is the induced map in Definition~\ref{dfn_poly prod functoriality}. Further, if the maps $f_i$ have the additional property of being homotopy equivalences of pairs (or homeomorphisms of pairs) for $1\leq i\leq m$, then so do the maps $f(\sigma,c)$ for every $\sigma\in\Delta^{m-1}$.
\end{proposition}

We call $f(\bullet,c)=\{f(\sigma,c)|\sigma\in\Delta^{m-1}\}$ \emph{the induced collection of maps}. Given a simplicial complex $K$ on $[m]$, define $f(K,c):(\underline{X},\underline{A})^{K,c}\to(\underline{Y},\underline{B})^{K,c}$ to be the colimit of
\begin{equation}\label{dfn_varphi K}
(\underline{X},\underline{A})^{\sigma,c}\overset{f(\sigma,c)}{\longrightarrow}(\underline{Y},\underline{B})^{\sigma,c}\overset{\widetilde{\imath}^{\,K_{\sigma}}}{\longrightarrow}(\underline{Y},\underline{B})^{K,c}
\end{equation}
over $\sigma\subsetneq[m]$. This is well defined due to Diagram~\ref{diagram_unique weight poly prod1}. One can use Diagram~\ref{diagram_unique weight poly prod2} to show there is a commutative diagram
\[
\xymatrix{
(\underline{X},\underline{A})^{K}\ar[r]^-{f(K)}\ar[d]^-{\eta(K)}	&(\underline{Y},\underline{B})^{K}\ar[d]^-{\widetilde{\eta}(\partial\sigma)}\\
(\underline{X},\underline{A})^{K,c}\ar[r]^-{f(K,c)}		&(\underline{Y},\underline{B})^{K,c}.
}
\]
%\stanley{took out suspension splitting since it is not used}

%\so{here I add the retraction property of weighted polyhedral products, which is used in Lemma~\ref{lemma_homology}} 

There are also systematic retractions off $(\underline{X},\underline{A})^{K,c}$.

\begin{proposition}[Corollary 5.6 in~\cite{SST}]\label{prop_weigted poly prod retract}
Let $\sigma$ be a simplex of $K\subset[m]$, and let $(\underline{X},\underline{A})^{\sigma,c}$ be the associated weighted polyhedral product without ghost vertices. Then the inclusion $\imath^K_{\sigma}\colon (\underline{X},\underline{A})^{\sigma,c}\to (\underline{X},\underline{A})^{K,c}$ has a left homotopy inverse 
$p^{\sigma}_K\colon (\underline{X},\underline{A})^{K,c}\to(\underline{X},\underline{A})^{\sigma,c}$.
\end{proposition}

\section{Setting up the realization problem for $A(\mathfrak{c},\underline{d})$ when $m=3$}
Given $m\in\N$, let $\underline{d}=\{d_1,\ldots,d_m\}$ be a sequence of positive integers, let $\mathfrak{c}=\{\mathfrak{c}_{\sigma}\mid\sigma\subseteq[m]\}$ be a coefficient sequence, and let $A(\mathfrak{c},\underline{d})$ be the graded-commutative ring defined by Definition~\ref{def:mgwspr} and Equation~\eqref{dfn_A(c) prod}. We wish to construct a realization of $A(\mathfrak{c},\underline{d})$. Note that in the case $\mathfrak{c}=\Phi(c)$ where $c$ is a power sequence this has already been done in \cite{SST}, as stated in Theorem~\ref{realization}. 

Let us start with the special case where $\mathfrak{c}_{\sigma}=1$ for all $\sigma\subseteq[m]$. Then 
by Equation~\eqref{dfn_A(c) prod}, $$A(\mathfrak{c},\underline{d})=\frac{F_\Z[x_1, \dots, x_m]}{(x_1^2, \dots, x_m^2)}$$
is the sphere-product algebra on the set  of generators $\{x_1,\ldots,x_m\}$ with $|x_i|=d_i$. It is geometrically realized 
by $S^{d_1}\times\cdots\times S^{d_m}$. Recall that $d_{[m]}=d_{1}+\cdots + d_{m}$. 
The $(d_{[m]}-1)$-skeleton of \linebreak $S^{d_1}\times\cdots\times S^{d_m}$ is the fat wedge $FW(S^{d_1},\ldots,S^{d_m})$ and its $d_{[m]}$-cell is attached to $FW(S^{d_1},\ldots,S^{d_m})$ by a higher Whitehead product. In Example~\ref{ex_higher wh prod} we see that, in the terms of polyhedral products, the boundary of {the top cell} $e^{d_{[m]}}$ and the fat wedge can be expressed as $(\underline{CS},\underline{S})^{\partial\Delta^{m-1}}$ and $(\underline{\Sigma S},\ast)^{\partial\Delta^{m-1}}$ respectively, while the attaching map is the induced map
\[
q(\partial\Delta^{m-1}):(\underline{CS},\underline{S})^{\partial\Delta^{m-1}}\to(\underline{\Sigma S},\ast)^{\partial\Delta^{m-1}}.
\]

When $m=3$ and it is not the case that $\cf_{\sigma}= 1$ for every $\sigma\subseteq [m]$, the strategy for proving Theorem~\ref{thm_top_main} is to modify the construction of $S^{d_1}\times S^{d_2}\times S^{d_3}$ by replacing the fat wedge by a weighted polyhedral product and the higher Whitehead product $q(\partial\Delta^2)$ by a suitable map.

Given a coefficient sequence $\mathfrak{c}=\{\mathfrak{c}_{\sigma}\in\N\mid\sigma\subseteq[3]\}$ and a sequence $\underline{d}=\{d_1,d_2,d_3\in\N\mid d_i>1\}$ of degrees, we construct a power sequence
\[
c:\CAT(\partial\Delta^2)\to\N^3,\quad c(\sigma)=(c^{\sigma}_1,c^{\sigma}_2, c^{\sigma}_3)
\]
by taking
\begin{equation}\label{eqn_power seq c}
c^{\sigma}_i=\begin{cases}
\mathfrak{c}_{\{1,2\}}	&\text{if }(\sigma,i)=(\{1,2\},1)\text{ or }(\sigma,i)=(\{1,2,3\},1)\\
\mathfrak{c}_{\{2,3\}}	&\text{if }(\sigma,i)=(\{2,3\},2)\text{ or }(\sigma,i)=(\{1,2,3\},2)\\
\mathfrak{c}_{\{1,3\}}	&\text{if }(\sigma,i)=(\{1,3\},3)\text{ or }(\sigma,i)=(\{1,2,3\},3)\\
1		&\text{other cases}.
\end{cases}
\end{equation}
Let $(\underline{CS},\underline{S})=\{(CS^{d_i-1},S^{d_i-1})\}^3_{i=1}$ and $(\underline{\Sigma S},\ast)=\{(S^{d_i},\ast)\}^3_{i=1}$ be sequences of $CW$- power couples whose power maps are described in Example~\ref{ex_cone}, and let
\[
\{(\underline{CS},\underline{S})^{\bullet,c},\eta(\bullet)\}\quad
\text{and}\quad
\{(\underline{\Sigma S},\ast)^{\bullet,c},\tilde{\eta}(\bullet)\}
\]
be their associated weighted polyhedral product systems. For $1\leq i\leq 3$ let
\[
q_i:(CS^{d_i-1},S^{d_i-1})\to(S^{d_i},\ast)
\]
be the quotient map collapsing the boundary $S^{d_i-1}$ of $CS^{d_i-1}$ to a point. Then each $q_i$ is a map between $CW$--power couples. By Proposition~\ref{lemma_map between weight poly prod} they induce a collection of maps
\begin{equation}\label{eqn_weighted poly prod quotient}
q(\bullet,c)=\{q(\sigma,c):(\underline{CS},\underline{S})^{\sigma,c}\to(\underline{\Sigma S},\ast)^{\sigma,c}|\sigma\subseteq[3]\}    
\end{equation}
making the diagrams
\begin{enumerate}[label=(\alph*)]
\begin{minipage}{0.45\linewidth}
\item\label{diagram_q incl}
\hspace{0.5cm}
$\xymatrix{
(\underline{CS},\underline{S})^{\tau,c}\ar[r]^-{q(\tau,c)}\ar[d]^-{\imath^{\sigma}_{\tau}}	&(\underline{\Sigma S},\ast)^{\tau,c}\ar[d]^-{\widetilde{\imath}^{\,\sigma}_{\tau}}\\
(\underline{CS},\underline{S})^{\sigma,c}\ar[r]^-{q(\sigma,c)}		&(\underline{\Sigma S},\ast)^{\sigma,c}
}$
\end{minipage}
\hfill
\begin{minipage}{0.45\linewidth}
\item\label{diagram_q n eta}
\hspace{0.5cm}
$\xymatrix{
(\underline{CS},\underline{S})^{\sigma}\ar[r]^-{q(\sigma)}\ar[d]^-{\eta(\sigma)}	&(\underline{\Sigma S},\ast)^{\sigma}\ar[d]^-{\widetilde{\eta}(\sigma)}\\
(\underline{CS},\underline{S})^{\sigma,c}\ar[r]^-{q(\sigma,c)}		&(\underline{\Sigma S},\ast)^{\sigma,c}
}$
\end{minipage}
\end{enumerate}
commute, where $q(\sigma)$ is the restriction of $\prod^3_{i=1}q_i$ to $(\underline{CS},\underline{S})^{\sigma}$. The colimit
\[
q(\partial\Delta^2,c):(\underline{CS},\underline{S})^{\partial\Delta^2,c}\longrightarrow(\underline{\Sigma S},\ast)^{\partial\Delta^2,c}
\]
makes the diagram
\[
\xymatrix{
(\underline{CS},\underline{S})^{\partial\Delta^2}\ar[rr]^-{q(\partial\Delta^2)}\ar[d]^-{\eta(\partial\Delta^2)}	&&(\underline{\Sigma S},\ast)^{\partial\Delta^2}\ar[d]^-{\widetilde{\eta}(\partial\Delta^2)}\\
(\underline{CS},\underline{S})^{\partial\Delta^2,c}\ar[rr]^-{q(\partial\Delta^2,c)}		&&(\underline{\Sigma S},\ast)^{\partial\Delta^2,c}
}
\]
commute. 

We want to compare the weighted polyhedral products $(\underline{CS},\underline{S})^{\partial\Delta^2,c}$ and $(\underline{\Sigma S},\ast)^{\partial\Delta^2,c}$ to the normal polyhedral products $(\underline{CS},\underline{S})^{\partial\Delta^2}$ and $(\underline{\Sigma S},\ast)^{\partial\Delta^2}$.
%\stanley{I took out (what was) Lemma 3.1 since it is not needed in the new proof}
To do so we determine the homotopy type of $(\underline{CS},\underline{S})^{\partial\Delta^2,c}$ in Proposition~\ref{prop_retraction splitting}. This first requires a homology calculation.

\begin{lemma}\label{lemma_homology}
The homology of $(\underline{CS},\underline{S})^{\partial\Delta^2,c}$ is given by
\[
H_j((\underline{CS},\underline{S})^{\partial\Delta^2,c};\Z)\cong\begin{cases}
\Z				&\text{if }j=0,d_{\{1,2,3\}}-1\\
\bigoplus_{s=1}^{t_{j}}\Z/ \mathfrak{c}_{\{i_1,i_2\}}\Z	&\text{if } j=d_{\{i_1,i_2\}}-1\\
\text{torsion}	&\text{if }j=d_{\{1,2,3\}}-2\\
0		&\text{else},
\end{cases}
\]
where for each $\{i_1,i_2\}\subset [3]$ we have $1\leq d_{\{i_1,i_2\}}-1\leq d_{\{1,2,3\}}-3$ and $t_{j}$ records the number of distinct pairs $\{i_1,i_2\}$ contributing nonzero homology in the same degree.   
Moreover, if $\ell=\text{lcm}( \mathfrak{c}_{\{1,2\}}, \mathfrak{c}_{\{1,3\}}, \mathfrak{c}_{\{2,3\}})$ then 
there are generators $u\in H_{d_{\{1,2,3\}}-1}((\underline{CS},\underline{S})^{\partial\Delta^2};\Z)\cong \Z$ and 
$v \in H_{d_{\{1,2,3\}}-1}((\underline{CS},\underline{S})^{\partial\Delta^2, c};\Z)$ such that 
$$\eta(\partial\Delta^2)(u)=\frac{ \mathfrak{c}_{\{1,2\}} \mathfrak{c}_{\{2,3\}} \mathfrak{c}_{\{1,3\}}}{\ell}v.$$
\end{lemma}

\begin{proof}

Let $K\subset\partial\Delta^2$ be a simplicial complex. We describe cell structures on $(\underline{CS},\underline{S})^K$ 
and $(\underline{CS},\underline{S})^{K,c}$ and compute their cellular chain complexes $C^\cell_*$. 
The cell structure for $(CS^{d_i-1}, S^{d_i})$ is the standard one with one $0$-cell, one $(d_i-1)$-cell and one
$d_i$-cell, giving the chain complex:
\[
C^\cell_*(CS^{d_i-1})=\Z\langle 1, x_i, y_i\rangle, \ dy_i=x_i, dx_i=d_1=0, \
|y_i|=d_i=|x_i|+1.
\]
The cell structure on {$CS^{d_1-1}\times CS^{d_2-1}\times CS^{d_3-1}$} is the product cell structure with a basis for cellular chains given by elements of the form $a_1\otimes a_2\otimes a_3$ where $a_i=1, x_i \mbox{ or } y_i$ and the differential is computed using the Koszul convention. 
We will sometimes not write $a_i$ if it is $1$.
For any $K$, 
$C_*^{\cell}((\underline{CS},\underline{S})^K)$ is the corresponding subcomplex.  

Let  $L=\{ \emptyset, \{1\}, \{2\}, \{3\}\}$. Then since $c_i^{\{ i\}}=1$ we can use the same the cell structure on 
$(\underline{CS},\underline{S})^{L,c}$ and consider 
$C_*^{\cell}((\underline{CS},\underline{S})^{L,c})=C_*^{\cell}((\underline{CS},\underline{S})^L)$. 
More explicitly, $C_*^{\cell}((\underline{CS},\underline{S})^L)$ has basis  

\[
\left\lbrace\begin{array}{c}
1,\,x_1,\,x_2,\,x_3,\,x_1\otimes x_2,\,x_1\otimes x_3,\,x_2\otimes x_3,x_1\otimes x_2\otimes x_3,\\
y_1,\,y_1\otimes x_2,\,y_1\otimes x_3,\,y_1\otimes x_2\otimes x_3,\\
y_2,\,x_1\otimes y_2,\,y_2\otimes x_3,\,x_1\otimes y_2\otimes x_3,\\
y_3,\,x_1\otimes y_3,\,x_2\otimes y_3,x_1\otimes x_2\otimes y_3
\end{array}\right\rbrace.
\]

This also works for subsimplicial complexes of $L$, so
$C_*^{\cell}((\underline{CS},\underline{S})^{\partial\{ 1,2\}})=C_*^{\cell}((\underline{CS},\underline{S})^{\partial\{ 1,2\},c})$ has a basis given by removing the terms with $y_3$ in them 
(so all except the last row). Observe that 
$$C_*^{\cell}((\underline{CS},\underline{S})^{\{ 1,2\}})=C_*^{\cell}((\underline{CS},\underline{S})^{\partial\{ 1,2\}})
\oplus \Z\langle y_1\otimes y_2, y_1\otimes y_2\otimes x_3\rangle$$ and there is a pushout
\[
\xymatrix{
  (\underline{CS},\underline{S})^{\partial \{ 1,2\}} \ar[r] \ar[d]_{\eta(\partial{\{ 1,2\}})} 
  & (\underline{CS},\underline{S})^{\{ 1,2\}} \ar[d]^{\eta( \{ 1,2\})} \\
  (\underline{CS},\underline{S})^{\partial \{ 1,2\}, c} \ar[r] & (\underline{CS},\underline{S})^{\{ 1,2\}, c}
}
\]
{
Observe that the top horizontal map is the canonical inclusion $S^{d_{\{1,2\}}-1}\times S^{d_3-1}\hookrightarrow D^{d_{\{1,2\}}}\times S^{d_3-1}$ and the left vertical map is
\[
\eta(\overline{\partial\{1,2\}})\times id\colon
S^{d_{\{1,2\}}-1}\times S^{d_3-1}\cong(\underline{CS},\underline{S})^{\overline{\partial\{1,2\}}}\times S^{d_3-1}\longrightarrow
(\underline{CS},\underline{S})^{\overline{\partial\{1,2\}},c}\times S^{d_3-1}
\]
where $\overline{\partial\{1,2\}}$ denotes the two disjoint points on the vertex set $\{1,2\}$.
Hence there is a pushout}
%Lemma \ref{lemma_colim commute po} implies there is a pushout 
\[
\xymatrix{
  C_*^{\cell}((\underline{CS},\underline{S})^{\partial\{ 1,2\}}) \ar[r] \ar[d]_{\eta(\partial \{ 1,2\})_*} 
  &  C_*^{\cell}((\underline{CS},\underline{S})^{\{ 1,2\}}) \ar[d]^{\eta(\{ 1,2\})_*} \\
   C_*^{\cell}((\underline{CS},\underline{S})^{\partial\{ 1,2\}, c}) \ar[r] &  C_*^{\cell}((\underline{CS},\underline{S})^{{\{ 1,2\}}, c})
}
\]
To compute $\eta(\partial {\{ 1,2\}})_*$ use Definition \ref{dfn_weighted poly prod system} (iii)(c). 
Note that $\sigma=\{ 1,2\}$, so for 
$\tau=\emptyset, \{1\} \mbox{ or } \{2\}$, we obtain $C_*^{\cell}(c^{\sigma/\tau})= \mathfrak{c}_{\{1,2\}}\otimes 1\otimes 1$. Therefore 
\[
\eta(\partial \{ 1,2\})_*(a_1\otimes a_2\otimes a_3)=
\begin{cases}
   \mathfrak{c}_{\{1,2\}}  a_1\otimes a_2\otimes a_3 & \text{if } a_1=x_1 \mbox{ or } y_1 \\
   a_2\otimes a_3   & \text{if } a_1=1.
\end{cases}
\]
Let $z_{12}$ and $z_{12}\otimes x_3$ denote the images of $y_1\otimes y_2$ and $y_1\otimes y_2 \otimes x_3$.
Then 
$$C_*^{\cell}((\underline{CS},\underline{S})^{\{ 1,2\}, c})=
C_*^{\cell}((\underline{CS},\underline{S})^{\partial\{ 1,2\}, c})\oplus \Z\langle z_{12}, z_{12}\otimes x_3\rangle$$ 
with $dz_{12}= \mathfrak{c}_{\{1,2\}}(x_1\otimes y_2+(-1)^{d_1}y_1\otimes x_2)$ and 
$d(z_{12}\otimes x_3)= \mathfrak{c}_{\{1,2\}}(x_1\otimes y_2\otimes x_3+(-1)^{d_1}y_1\otimes x_2\otimes x_3)$. Also 
\[
\eta(\{ 1,2\})_*(\alpha)=
\begin{cases}
  \eta(\partial\{ 1,2\})_*(\alpha) & \text{if } \alpha\in  C_*^{\cell}((\underline{CS},\underline{S})^{\partial\{ 1,2\}})\\
   z_{12}  & \text{if } \alpha=y_1\otimes y_2\\
   z_{12}\otimes x_3  & \text{if } \alpha=y_1\otimes y_2\otimes x_3
\end{cases}
\]
Similarly,
%\stanley{should we add the commutative diagrams corresponding to adding two more simplices here, or is it already clear?}
%\so{I think the cellular chains are quite clear and no need to add extra commutative diagrams}
we compute $$C_*^{\cell}((\underline{CS},\underline{S})^{\partial \{ 1,2,3\}, c})
=C_*^{\cell}((\underline{CS},\underline{S})^{L, c})\oplus \Z\langle z_{12}, z_{12}\otimes x_3, 
z_{23}, x_1\otimes z_{23}, z_{13}, z_{13}\widetilde{\otimes} x_2\rangle$$ 
where: $z_{23}$ and $x_1\otimes z_{23}$ are the images under $\eta(\{ 2,3\})_*$ (followed by the inclusion) of $y_2\otimes y_3$ and $x_1\otimes y_2\otimes y_3$; similarly for $z_{13}$ and $z_{13}\widetilde{\otimes} x_2$, with  $\widetilde{\otimes}$ referring to the tensor product with entries $2$ and~$3$ swapped; and the differential is given by
\[\begin{split} 
    dz_{23} & =\mathfrak{c}_{\{2,3\}}(x_2\otimes y_3 +(-1)^{d_2} y_2\otimes x_3) \\  
    d(x_1\otimes z_{23}) & =\mathfrak{c}_{\{2,3\}}((-1)^{d_1-1}x_1\otimes x_2\otimes y_3+(-1)^{d_1+d_2-1}x_1\otimes y_2\otimes x_3) \\ 
    dz_{13} & =\mathfrak{c}_{\{1,3\}}(x_1\otimes y_3 +(-1)^{d_1} y_1\otimes x_3) \\  
    d(z_{13}\widetilde{\otimes} x_2) & =\mathfrak{c}_{\{1,3\}}(x_1\otimes x_2\otimes y_3+(-1)^{d_1+d_2-1}y_1\otimes x_2\otimes x_3). 
\end{split}\] 
A calculation similar to that for $\eta(\partial\{1,2\})$ gives 
\[
\eta(\partial\{ 1,2,3\})_*(y_1\otimes y_2\otimes x_3)= \mathfrak{c}_{\{1,3\}} \mathfrak{c}_{\{2,3\}}z_{12}\otimes x_3, 
\]
\[
\eta(\partial\{ 1,2,3\})_*(y_1\otimes x_2\otimes y_3)= \mathfrak{c}_{\{1,2\}} \mathfrak{c}_{\{2,3\}}z_{13}\widetilde{\otimes} x_2,
\]
\[
\eta(\partial\{ 1,2,3\})_*(x_1\otimes y_2\otimes y_3)= \mathfrak{c}_{\{1,2\}} \mathfrak{c}_{\{1,3\}}x_1\otimes z_{23}.
\]
Having an explicit chain complex, it is straightforward to compute that $H_*((\underline{CS},\underline{S})^{\partial \{ 1,2,3\}, c};\Z)$ is as stated in the lemma. 

Further, if $\ell=\text{lcm}( \mathfrak{c}_{\{1,2\}}, \mathfrak{c}_{\{1,3\}}, \mathfrak{c}_{\{2,3\}})$, let 
\[\begin{split} 
  u & =(-1)^{d_1+d_2}y_1\otimes y_2\otimes x_3+(-1)^{d_1}y_1\otimes x_2\otimes y_3+x_1\otimes y_2\otimes y_3 \\  
  v & =(-1)^{d_1+d_2}\frac{\ell}{ \mathfrak{c}_{\{1,2\}}}z_{12}\otimes x_3+(-1)^{d_1}\frac{\ell}{ \mathfrak{c}_{\{1,3\}}}z_{13}\widetilde{\otimes} x_2+\frac{\ell}{ \mathfrak{c}_{\{2,3\}}}x_1\otimes z_{23}. 
\end{split}\] 
Observe that $\eta(\partial\{ 1,2,3\})_*(u)=\frac{ \mathfrak{c}_{\{1,2\}} \mathfrak{c}_{\{2,3\}} \mathfrak{c}_{\{1,3\}}}{\ell}v$ and that $u$ and $v$ generate $H_{d_{\{1,2,3\}}-1}((\underline{CS},\underline{S})^{\partial\Delta^2};\Z)$ and $H_{d_{\{1,2,3\}}-1}((\underline{CS},\underline{S})^{\partial\Delta^2, c};\Z)$. 
\end{proof}

The description of the homology of 
$(\underline{CS},\underline{S})^{\partial\Delta^2,c}$ 
in Lemma~\ref{lemma_homology} lets us determine its homotopy type.

\begin{proposition}\label{prop_retraction splitting}
Let $G=H_{d_{\{1,2,3\}}-2}((\underline{CS},\underline{S})^{\partial\Delta^2,c};\Z)$.
If $d_1,d_2,d_3\geq2$ then there is a homotopy equivalence
\[
(\underline{CS},\underline{S})^{\partial\Delta^2,c}\simeq S^{d_{\{1,2,3\}}-1}\vee P^{d_{\{1,2,3\}}-1}(G)\vee\bigvee_{\{i_1,i_2\}\subset[3]}\left(P^{d_{\{i_1,i_2\}}}( \mathfrak{c}_{\{i_1,i_2\}})\right).
\]
\end{proposition}

\begin{proof}
Let $i_{1},i_{2}\in\{1,2,3\}$ with $1\leq i_{1}<i_{2}\leq 3$. Note that there are three possibilities: 
$(1,2)$, $(1,3)$ and $(2,3)$. Since the face $(i_{1},i_{2})$ is a 
full subcomplex of $\partial\Delta^{2}$, by Proposition~\ref{prop_weigted poly prod retract} the inclusion 
\(\namedright{(\underline{CS},\underline{S})^{(i_{1},i_{2}),c}}{}{(\underline{CS},\underline{S})^{\partial\Delta^2,c}}\) 
has a left homotopy inverse. Given two distinct faces $(i_{1},i_{2})$ and $(i'_{1},i'_{2})$ 
these will intersect at one vertex, say $i_{1}=i'_{1}$ (the other cases being similar). Since the left 
homotopy inverses in Proposition~\ref{prop_weigted poly prod retract} are induced by projection maps, the retractions of 
$(\underline{CS},\underline{S})^{(i_{1},i_{2}),c}$ and $(\underline{CS},\underline{S})^{(i'_{1},i'_{2}),c}$ 
off $(\underline{CS},\underline{S})^{\partial\Delta^2,c}$ will intersect at $(\underline{CS},\underline{S})^{\{i_{1}\},c}$. 
By assumption, $c_{i_{1}}=1$, so by Example~\ref{eg_ctrivial}, 
$(\underline{CS},\underline{S})^{\{i_{1}\},c}=(\underline{CS},\underline{S})^{\{i_{1}\}}=CS^{d_{i_{1}}}$, 
which is contractible. Thus, up to homotopy, the retractions of 
$(\underline{CS},\underline{S})^{(i_{1},i_{2}),c}$ and $(\underline{CS},\underline{S})^{(i'_{1},i'_{2}),c}$ 
off $(\underline{CS},\underline{S})^{\partial\Delta^2,c}$ meet only at the basepoint. 

By Example~\ref{eg_Moore sp}, there is a homotopy equivalence 
$(\underline{CS},\underline{S})^{(i_{1},i_{2}),c}\simeq P^{d_{\{i_{1},i_{2}\}}}( \mathfrak{c}_{\{i_{1},i_{2}\}})$. 
We may take the wedge sum of the three inclusion maps 
\(\namedright{P^{d_{\{i_1,i_2\}}}( \mathfrak{c}_{\{i_1,i_2\}})}{}{(\underline{CS},\underline{S})^{\partial\Delta^{2},c}}\) 
to obtain 
\[a\colon\namedright{\bigvee_{\{i_{1},i_{2}\}\subset [3]} P^{d_{\{i_1,i_2\}}}(c_{\{i_1,i_2\}})}{} 
     {(\underline{CS},\underline{S})^{\partial\Delta^{2},c}}.\] 
As $(\underline{CS},\underline{S})^{\partial\Delta^{2},c}$ need not be a co-$H$-space, it may not be possible to add the 
three left homotopy inverses. Instead, taking their product, there is 
a homotopy commutative diagram 
\begin{equation} 
  \label{retract1} 
  \diagram 
      \bigvee_{\{i_{1},i_{2}\}\subset [3]} P^{d_{\{i_1,i_2\}}}( \mathfrak{c}_{\{i_1,i_2\}})\rto^-{a}\drto^{\iota} 
           & (\underline{CS},\underline{S})^{\partial\Delta,c}\dto^{p} \\ 
       & \prod_{\{i_{1},i_{2}\}\subset [3]} P^{d_{\{i_1,i_2\}}}( \mathfrak{c}_{\{i_1,i_2\}}) 
  \enddiagram 
\end{equation}  
where $p$ is the product of the three left homotopy inverses and $\iota$ is the inclusion of the wedge into the product. 

Consider the map $p$. By Lemma~\ref{lemma_homology}, the dimension of 
$(\underline{CS},\underline{S})^{\partial\Delta^{2},c}$ 
is $d_{\{1,2,3\}}-1=d_{1}+d_{2}+d_{3}-1$. 
On the other hand, let $C_\iota$ be the homotopy cofibre 
of $\iota$. As $\iota$ is the inclusion of a wedge into the product, 
the connectivity of $C_\iota$ equals the connectivity of the smash 
product of the two Moore space factors in 
$\prod_{\{i_1,i_2\}\subset [3]} P^{d_{\{i_1,i_2\}}}( \mathfrak{c}_{\{i_1,i_2\}})$ 
of least connectivity. Thus $C_\iota$ has connectivity equal to the minimum 
of the numbers $d_{\{i_{1},i_{2}\}}+d_{\{i'_{1}i'_{2}\}}-3$,
where $\{i_{1},i_{2}\}$ and $\{i'_{1},i'_{2}\}$ are distinct pairs in $[3]$. 
Any such distinct pairs has union $\{1,2,3\}$ and intersect at one element. Thus 
$d_{\{i_{1},i_{2}\}}+d_{\{i'_{1},i'_{2}\}}-3\geq d_{i_{1}}+d_{i_{2}}+d_{i_{3}}+d_{j}-3$, 
where $d_{j}$ is the minimum of $d_{i_{1}}$, $d_{i_{2}}$ and $d_{i_{3}}$. Since 
each $d_{i_{k}}\geq 2$, the connectivity of $C_\iota$ is at least as great 
as the dimension of $(\underline{CS},\underline{S})^{\partial\Delta^{2},c}$. Thus there 
is a lift 
\begin{equation} 
  \label{retract2} 
  \diagram 
      & (\underline{CS},\underline{S})^{\partial\Delta^{2},c}\dto^{p}\dlto_{b} \\ 
      \bigvee_{\{i_{1},i_{2}\}\subset [3]} P^{d_{\{i_{1},i_{2}\}}}( \mathfrak{c}_{\{i_1,i_2\}})\rto^-{\iota} 
           & \prod_{\{i_{1},i_{2}\}\subset [3]} P^{d_{\{i_1,i_2\}}}( \mathfrak{c}_{\{i_1,i_2\}}) 
  \enddiagram 
\end{equation}  
for some map $b$. The homotopy commutativity of~(\ref{retract1}) and~(\ref{retract2}) 
imply that the composite $b\circ a$ induces an isomorphism in homology and so is a 
homotopy equivalence.

The description of $H_{\ast}((\underline{CS},\underline{S})^{\partial\Delta^{2},c};\Z)$ in 
Lemma~\ref{lemma_homology} implies that there is a homotopy cofibration  
\[\nameddright{S^{d_{\{1,2,3\}}-2}\vee P^{d_{\{1,2,3\}}-2}(G)}{f} 
     {\bigvee_{\{i_{1},i_{2}\}\subset [3]} P^{d_{\{i_1,i_2\}}}( \mathfrak{c}_{\{i_1,i_2\}})} 
     {a}{(\underline{CS},\underline{S})^{\partial\Delta^{2},c}}\] 
where $f$ attaches the top cells to $(\underline{CS},\underline{S})^{\partial\Delta^{2},c}$. 
Since $a$ has a left homotopy inverse, $f$ is null homotopic. Therefore there is a homotopy 
equivalence 
\[
(\underline{CS},\underline{S})^{\partial\Delta^2,c}\simeq S^{d_{\{1,2,3\}}-1}\vee P^{d_{\{1,2,3\}}-1}(G)\vee\bigvee_{\{i_1,i_2\}\subset[3]}P^{d_{\{i_1,i_2\}}}( \mathfrak{c}_{\{i_1,i_2\}}).
\]
\end{proof}

\section{The proof of Theorem~\ref{thm_top_main}}

In~\cite[Theorem 8.4]{SST} it was proved that an $m$-generated weighted sphere product ring $A(\mathfrak{c},\underline{d})$ can be realized as the cohomology of a weighted polyhedral product $(\underline{\Sigma S},\ast)^{\Delta^{m-1},c}$ for the subset of coefficient sequences $\mathfrak{c}$ that are power sequences. Now when $m=3$ we want to realize $A(\mathfrak{c},\underline{d})$ for the larger class of coefficient sequences described in Equation (\ref{eqn_power seq c}).
The idea is to use $(\underline{\Sigma S},\ast)^{\Delta^2,c}$ to construct a realization.

\begin{lemma}\label{lemma_cohmgly ring of (Sigma S,*)}

Let $(\underline{\Sigma S},\ast)=\{(S^{d_i},\ast)\}^3_{i=1}$ and let $c$ be defined in~\eqref{eqn_power seq c}.
Then as graded modules
\[
H^*((\underline{\Sigma S},\ast)^{\partial\Delta^2,c};\Z)\cong\bigoplus_{\sigma\subsetneq[3]}\langle{a_{\sigma}}\rangle
\quad\text{and}\quad
H^*((\underline{\Sigma S},\ast)^{\Delta^2,c};\Z)\cong\bigoplus_{\sigma\subseteq[3]}\langle{a_{\sigma}}\rangle,
\]
where $a_{\sigma}$ has degree $d_{\sigma}$, and the multiplication is given by
\[
a_ia_j=\mathfrak{c}_{\{i,j\}}a_{\{i,j\}},\quad
a_1a_2a_3=\mathfrak{c}_{\{1,2\}}\mathfrak{c}_{\{2,3\}}\mathfrak{c}_{\{1,3\}}a_{\{1,2,3\}}.
\]
\end{lemma}

\begin{proof}
This follows directly from~\cite[Theorems~7.12 and~8.4]{SST} and the definition of $c$ in~\eqref{eqn_power seq c}.  
\end{proof}

\begin{lemma}\label{gooddiag}
Let $q(\bullet)\colon(\underline{CS},\underline{\Sigma S})^{\bullet}\to(\underline{\Sigma S},\ast)^{\bullet}$ and $q(\bullet,c)\colon(\underline{CS},\underline{S})^{\bullet,c}\to(\underline{S},\ast)^{\bullet,c}$ be induced by the collapsing maps $CS^{d_i-1}\to S^{d_i}$ given in Example~\ref{ex_higher wh prod} and~\eqref{eqn_weighted poly prod quotient}.
Consider the diagram

\[
\xymatrix{
  &&& S^{d_{\{1,2,3\}}-1}\cong (\underline{CS},\underline{S})^{\partial\Delta^2} 
  \ar[rr] \ar[d]^-{q(\partial\Delta^2)} \ar[llldd]_{\eta(\partial \Delta^2)} &&  (\underline{CS},\underline{S})^{\Delta^2}  \ar[d]^-{q(\Delta^2)} \\
  &&& (\underline{\Sigma S},*)^{\partial\Delta^2} \ar[rr] \ar[d]^{\eta(\partial \Delta^2)} 
  && (\underline{\Sigma S},*)^{\Delta^2} \ar[d]^{\eta(\Delta^2)} \\
  (\underline{CS},\underline{S})^{\partial\Delta^2,c} \ar[rrr]^{q(\partial\Delta^2,c)} &&& (\underline{\Sigma S},*)^{\partial\Delta^2,c}\ar[rr]
  && (\underline{\Sigma S},*)^{\Delta^2,c}
}
\]
The diagram commutes and the two right squares, and thus the outside right rectangle, are pushouts.

\end{lemma}
\begin{proof}
The upper right square is a pushout as explained in Example \ref{ex_higher wh prod} and the lower right square is a pushout by Definition \ref{dfn_weighted poly prod system} (iii). The commutativity of the triangle is given by 
Proposition \ref{lemma_map between weight poly prod}.
\end{proof}

By Example~\ref{ex_higher wh prod}, $S^{d_{\{1,2,3\}}-1}$ can be identified with $(\underline{CS},\underline{S})^{\partial\Delta^2}$. Let $u$ and $v$ 
be the generators from Lemma 
\ref{lemma_homology} and let $f\colon S^{d_{\{1,2,3\}}-1} \rightarrow 
(\underline{CS},\underline{S})^{\partial\Delta^2,c} $ be any map such that $f_*(u)=v$. Such a map exists because of the splitting in Proposition \ref{prop_retraction splitting}.

Let $\ell=\text{lcm}(\mathfrak{c}_{\{1,2\}}, \mathfrak{c}_{\{2,3\}}, \mathfrak{c}_{\{1,3\}})$. Then $\ell$ divides $\mathfrak{c}_{\{1,2,3\}}$ since, by definition of a coefficient sequence, each $\mathfrak{c}_{\{i,j\}}$ divides $\mathfrak{c}_{\{1,2,3\}}$ for all $i,j$.
Define $X$ as the cofiber of the composite 
$$S^{d_{\{1,2,3\}}-1}\stackrel{\frac{\mathfrak{c}_{\{1,2,3\}}}{\ell}}{\longrightarrow} 
S^{d_{\{1,2,3\}}-1}\stackrel{f}{\longrightarrow} 
(\underline{CS},\underline{S})^{\partial\Delta^2,c}\stackrel{q(\partial\Delta^2,c)}{\llarrow} (\underline{\Sigma S},*)^{\partial\Delta^2,c}.$$

\begin{proof}[Proof of Theorem~\ref{thm_top_main}] 
We will show that there is a ring isomorphism 
$H^*(X;\Z)\cong A(\mathfrak{c},\underline{d})$.
By Lemma \ref{lemma_homology}
 we know that $H_{d_{\{1,2,3\}}-1}(\eta(\partial \Delta^2))=\frac{\mathfrak{c}_{\{1,2\}}\mathfrak{c}_{\{2,3\}}\mathfrak{c}_{\{1,3\}}}{\ell}H_{d_{\{1,2,3\}}-1}(f)$.
Thus, since $$\pi_{d_{\{1,2,3\}}-1}((\underline{CS},\underline{S})^{\partial\Delta^2,c})\cong \Z \oplus \mbox{torsion}$$ 
by Proposition \ref{prop_retraction splitting}, there is a positive integer $k$ such that 
$\eta(\partial \Delta^2)\circ k\simeq f \circ k\frac{\mathfrak{c}_{\{1,2\}}\mathfrak{c}_{\{2,3\}}\mathfrak{c}_{\{1,3\}}}{\ell}$. 
Therefore there is a homotopy commutative diagram  
\[
\xymatrix{
 S^{d_{\{1,2,3\}}-1}  \ar[rrrr]^{q(\partial\Delta^2,c)\circ\eta(\partial \Delta^2)}  &&&& 
 (\underline{\Sigma S},*)^{\partial\Delta^2,c} \ar[r]  & (\underline{\Sigma S},*)^{\Delta^2,c}\\
S^{d_{\{1,2,3\}}-1}   \ar[rrrr]_{q(\partial\Delta^2,c)\circ f \circ k\frac{\mathfrak{c}_{\{1,2\}}\mathfrak{c}_{\{2,3\}}\mathfrak{c}_{\{1,3\}}\mathfrak{c}_{\{1,2,3\}}}{\ell^2}} \ar[u]^{k\frac{\mathfrak{c}_{\{1,2,3\}}}{\ell}} \ar[d]_{k\frac{\mathfrak{c}_{\{1,2\}}\mathfrak{c}_{\{2,3\}}\mathfrak{c}_{\{1,3\}}}{\ell}}
&&&& (\underline{\Sigma S},*)^{\partial\Delta^2,c} \ar[u]_= \ar[d]^= 
\ar[r] & X' \ar[d] \ar[u] \\ 
S^{d_{\{1,2,3\}}-1}  \ar[rrrr]_{q(\partial\Delta^2,c)\circ f \circ (\frac{\mathfrak{c}_{\{1,2,3\}}}{\ell})} 
&&&& (\underline{\Sigma S},*)^{\partial\Delta^2,c} \ar[r] & X\\ 
}
\]
where $X'$ is defined as the cofibre of $q(\partial\Delta^2,c)\circ f \circ k\frac{\mathfrak{c}_{\{1,2\}}\mathfrak{c}_{\{2,3\}}\mathfrak{c}_{\{1,3\}}\mathfrak{c}_{\{1,2,3\}}}{\ell^2}$ and the vertical maps in the third column are extensions that make the diagram commute up to homotopy. 
Lemma \ref{gooddiag} implies that the top row is a cofibration sequence, and the bottom row is a cofibration sequence by the definition of $X$. Note that $H^{i}((\underline{\Sigma S},\ast)^{\partial\Delta^2,c})=0$ for $i=d_{\{1,2,3\}}-1$ or $d_{\{1,2,3\}}$ by Lemma~\ref{lemma_cohmgly ring of (Sigma S,*)}.
% all the horizontal maps out of the sphere induce $0$ on ${\tilde{H}^*}$.
Thus we obtain a commutative diagram from the associated long exact sequences on cohomology:
\begin{equation}\label{homology_chase_diagram}
\xymatrix{
H^{d_{\{1,2,3\}}}(\underline{\Sigma S},*)^{\Delta^2,c}  \ar[d]
 & H^{d_{\{1,2,3\}}-1}(S^{d_{\{1,2,3\}}-1})\cong \Z \ar[d]^{k\frac{\mathfrak{c}_{\{1,2,3\}}}{\ell}} \ar[l]_{\cong} \\
H^{d_{\{1,2,3\}}}( X')  &H^{d_{\{1,2,3\}}-1}(S^{d_{\{1,2,3\}}-1}) \cong \Z   \ar[l]_{\cong} \\ 
H^{d_{\{1,2,3\}}}(X)\ar[u]& H^{d_{\{1,2,3\}}-1}(S^{d_{\{1,2,3\}}-1})\cong \Z   \ar[l]_{\cong}
\ar[u]_{k\frac{\mathfrak{c}_{\{1,2\}}\mathfrak{c}_{\{2,3\}}\mathfrak{c}_{\{1,3\}}}{\ell}}.\\ 
}
\end{equation} 
Since below dimension $d_{\{1,2,3\}}$ all of $(\underline{\Sigma S},*)^{\Delta^2,c}$, $X$ and $X'$ are equal to 
$(\underline{\Sigma S},*)^{\Delta^2,c}$ they all have cohomology linearly isomorphic to 
$\bigoplus_{\sigma\subseteq [m]} \Z\langle a_{\sigma}\rangle $. By
Lemma~\ref{lemma_cohmgly ring of (Sigma S,*)},
in all cases $a_ia_j=\mathfrak{c}_{\{i,j\}}a_{\{i,j\}}$, and in $H^*((\underline{\Sigma S},*)^{\Delta^2,c})$ we have  
$a_1a_2a_3=\mathfrak{c}_{\{1,2\}}\mathfrak{c}_{\{2,3\}}\mathfrak{c}_{\{1,3\}}a_{\{1,2,3\}}$. 

A diagram chase using the upper square in~(\ref{homology_chase_diagram}) implies that
$a_1a_2a_3=k\frac{\mathfrak{c}_{\{1,2,3\}}}{\ell}\mathfrak{c}_{\{1,2\}}\mathfrak{c}_{\{2,3\}}\mathfrak{c}_{\{1,3\}}a_{\{1,2,3\}}$ in $H^*(X')$, and another diagram chase using the lower square in~(\ref{homology_chase_diagram}) implies that  $a_1a_2a_3=\mathfrak{c}_{\{1,2,3\}}a_{\{1,2,3\}}$ in $H^*(X)$. Hence there is a ring isomorphism $H^*(X)\cong A(\mathfrak{c},\underline{d})$. Now take $X(\cf,\du)$ in the statement of Theorem~\ref{thm_top_main} to be $X$.
\end{proof}

\section{Orders in sphere product algebras with three generators}
\label{sec:algebras-with-three}

\subsection{Some facts about lattices}
\label{sec:some-facts-about}

If $V$ is a finite-dimensional $\Q$-vector space, recall that a subgroup $L \subset V$ with basis $\mathcal{B}$ is said to be a \emph{lattice} if $L$ is the smallest subgroup of $V$ containing $\mathcal{B}$. In the graded case, we require the basis $\mathcal B$ to consist of homogeneous elements of possibly varying degrees. If $M$ is a free abelian group of finite rank then $M$ is a lattice in $\Q \tensor_\Z M$, which is graded if $M$ is.

If $L$ is a lattice in $V$, then for every element $v \in V$, we may clear denominators to find a nonzero integer $d$ for which $dv \in L$.
Conversely, suppose $L$ is a subgroup of $V$ with the property that for every element $v \in V$, we may find a nonzero integer $d$ for which $dv \in L$. Then the canonical morphism $f: \Q \tensor_\Z L \to V$ is surjective. Note also that the rank of $L$ cannot exceed the dimension of $V$, so that $f$ is a surjective morphism from one $\Q$-vector space to another of greater or equal dimension. Therefore $f$ is an isomorphism and we conclude that $L$ is a lattice in $V$.\

\begin{lemma} \label{lem:torsionfree}
  Let $V$ be a $\Q$-vector space and let $W \subseteq U \subseteq V$ be vector subspaces. Let $M \subset V$ be a lattice in $V$. Then the quotient $(M \cap U)/(M\cap W)$ is torsion free.
\end{lemma}
\begin{proof}
  Let $x$ be an element of $M \cap U$ that has torsion image in $(M \cap U)/(M\cap W)$. That is, for some positive integer $d$, the element $dx$ lies in $M \cap W$. This implies $dx$ lies in $W$, which is a $\Q$-vector space, so $x$ lies in $W$, and consequently the image of $x$ in $(M \cap U)/(M\cap W)$ is $0$.
\end{proof}

\subsection{Structure of $A$}
\label{sec:structure-a}

% If $M$ is a free abelian group, then a \emph{basis} for $M$ is a set $\{b_i\}_{i \in I}$ for which $M = \bigoplus_{i \in I} \Z \langle b_i \rangle$

Fix a sequence of positive integers $\du = (d_1, d_2, d_3)$. Let $R$ denote the \mbox{$\Q$-algebra}
\[ R = \frac{F_\Q[x_1,x_2,x_3]}{(x_1^2, x_2^2, x_3^2)}, \quad \text{ where $|x_i|=d_i$.}\]
The relation $x_i^2 = 0$ is redundant if $d_i$ is odd. The algebra $R$ is connected and augmented. Let us write $N$ for the augmentation ideal, i.e., $(x_1, x_2, x_3)$.

For the sake of brevity, we will write $R_{123}$ instead of $R_{d_{\{1,2,3\}}}$ or $R_{d_1+d_2+d_3}$ for the degree-$(d_1+d_2+d_3)$ part of $R$, and similarly for other degrees.

Let $A$ be a graded order in $R$, as defined in the introduction. Our main structural result for $A$ is the following.
\begin{proposition} \label{pr:Astr}
  There exists a decomposition of finitely generated free abelian groups 
  \[  A = \Z \oplus L(1) \oplus L(2) \oplus L(3) \]
  in which
  \begin{enumerate}
  \item $N \cap A = L(1) \oplus L(2) \oplus L(3)$;
  \item $N^2 \cap A = L(2) \oplus L(3)$;
  \item $N^3 \cap A = L(3)$.
  \end{enumerate}
  Furthermore, if $\{y_1, y_2, y_3\}$ constitutes a basis for $L(1)$, then $y_1, y_2, y_3$ generate $R$ as a $\Q$-algebra.
\end{proposition}

\begin{proof}
  As part of the grading of $A$, there exists a decomposition $A = \Z \oplus (A \cap N)$, where the first summand consists of the degree-$0$ elements, and $(A \cap N)$ is generate by the positive-degree elements.  
  
  We next decompose $A \cap N$. Let us begin with $P$ defined by the short exact sequence
  \[
    \begin{tikzcd}
      0 \rar & N^2 \cap A \rar & N \cap A \rar & P \rar & 0 
    \end{tikzcd}
  \]
  where Lemma \ref{lem:torsionfree} assures us that $P$ is torsion free, and therefore free. The sequence admits a splitting in the category of graded abelian groups, and we choose such a splitting. Let us write $L(1)$ for the isomorphic image of $P$ in $N \cap A$. The decomposition of $N^2 \cap A$ into $L(2) \oplus L(3)$ is established similarly.

The last assertion of the proposition is proved as follows. If $\{y_1, y_2, y_3\}$ constitutes a basis for the abelian group $L(1)$, then the images of $y_1, y_2, y_3$ serve as a basis of the $\Q$-vector space $N/N^2$. This vector space also has the images of $x_1, x_2, x_3$ as a basis. In particular, for each $j \in \{1,2,3\}$, we can write some expression $x_j + p$ as a linear combination of the $y_i$, where $p$ is an expression in $N^2$. From here, it is easy to express each $x_i$ as a polynomial in $y_1, y_2, y_3$, which suffices since the $x_i$'s generate $R$ as a $\Q$-algebra.
\end{proof}

\begin{remark}
  For degree reasons, the abelian group $L(3)$ is exactly $A_{123}$. It is of rank $1$ and is generated by some nonzero rational multiple of $x_1x_2x_3$.
\end{remark}

\begin{remark}
  The ranks of the graded summands of $L(i)$ can be determined after applying $\Q \tensor_\Z -$, and since $\Q \tensor_\Z L(i) = N^i/N^{i+1}$, these ranks depend only on the $\Q$-algebra $R$ and are independent of $A$ and all choices made. In particular, a homogeneous basis for $L(1)$ must have the same degrees as $\{x_1, x_2, x_3\}$ (since the image of this set forms a basis for $N/N^2$), that is, $\du$ in some order.
\end{remark}

\subsection{Proof of \protect{Theorem \ref{th:ordersAllBasisG}}}\label{sec:proof-prot-refth}
 
\begin{proof}[Proof of \protect{Theorem \ref{th:ordersAllBasisG}}]
   Write down a decomposition $A = \Z \oplus L(1) \oplus L(2) \oplus L(3)$ as in Proposition~\ref{pr:Astr}.
  
   Endow $L(1)$ with an ordered homogeneous basis $\{y_1, y_2, y_3\}$, having degrees $(d_1, d_2,d_3)$ in this order. We claim that there exists a (possibly different) homogeneous basis $\{b_1, b_2, b_3\}$ of $L(1)$, satisfying $b_1^2=b_2^2=b_3^2=0$ and integers $c_{\{i,j\}}$ so that
   \[ \left\{ \frac{1}{c_{\{1,2\}}} b_1 b_2, \frac{1}{c_{\{2,3\}}}{b_2b_3}, \frac{1}{c_{\{1,3\}}}{b_1b_3}\right\}  \]
   constitutes a basis for $L(2)$.

  Multiplication gives us a morphism $L(1) \tensor_\Z L(1) \to N^2 \cap A$, whose image lies in degrees less than $d_{123}$, that is, a morphism
   \[ \mu : L(1) \tensor_\Z L(1) \to L(2). \]
   Observe that $ \mu$ becomes surjective after application of $\Q \tensor_\Z -$, since $\Q \tensor_\Z A$ is generated by $L(1)$ as an algebra.

   Let $M$ be the image of $\mu$, which is necessarily a rank-$3$ free abelian group. If $\mathcal B = \{b_1, b_2, b_3\}$ is a homogeneous basis of $L(1)$ that satisfies $b_1^2=b_2^2=b_3^2=0$, then $\{b_1b_2, b_2b_3, b_3b_1\}$ constitutes a generating set for $M$, which therefore constitutes an (ordered) basis $\mathcal B'$ of $M$, which we will say is the basis of $M$ \emph{induced by} $\mathcal B$. If we additionally choose an ordered basis $\mathcal Z$ of $L(2)$, then relative to~$\mathcal B'$ and $\mathcal Z$, the inclusion $M \subset L(2)$ is given by a $3 \times 3$ integer matrix $X$.

   Our claim may be restated as follows: there are choices of $\mathcal B$ as above and $\mathcal Z$ with respect to which $X$ takes the form
   \[
     \begin{bmatrix}
       c_{\{1,2\}} & 0 & 0 \\
       0 & c_{\{2,3\}} & 0 \\
       0 & 0 & c_{\{3,1\}}
     \end{bmatrix} \quad \text{ for some positive integers $c_{\{i,j\}}$.}\]
   We divide the proof of this into three cases.
   \begin{description}
   \item[The $d_i$ are all equal]
     Since the $d_i$  are equal, they are odd by hypothesis. Therefore the condition $b_1^2=b_2^2=b_3^2=0$ for a homogeneous basis will follow necessarily from graded-commutativity of the ambient ring, and does not require further verification. 
     
     Since $L(1)$ is an odd-degree summand of a graded-commutative algebra, we may identify $M = \Alt^2(L(1))$. For material on the alternating square $\Alt^2$, see Appendix~\ref{sec:altern-repr}.

     Changes of basis of $L(1)$, by $Y \in \GL(3; \Z)$, and of $L(2)$, by $Z \in \GL(3;\Z)$, change $X$ according to the rule
  \[ X \leadsto Z X \Alt^2(Y) \]
  and since $\Alt^2 : \GL(3;\Z) \to \SL(3;\Z)$ is surjective, according to Proposition \ref{pr:Zn=2r=3homo}, we may bring $X$ to  Smith normal form
  \begin{equation*}
    ZX \Alt^2(Y) =
    \begin{bmatrix}
      c_{\{1,2\}} & 0 & 0  \\ 0 & c_{\{2,3\}} & 0  \\ 0 & 0  & c_{\{1,3\}} 
    \end{bmatrix},
  \end{equation*}
 as required.
\item[Two of the $d_i$ are equal, the third is distinct] Without loss of generality, suppose that $d_1=d_2 \neq d_3$. Under our hypotheses, $d_1=d_2$ is odd. Therefore any homogeneous basis of $L(1)$ consists of two elements $b_1, b_2$ in odd degrees and $b_3 = \pm y_3$. The conditions $b_1^2=b_2^2=b_3^2=0$ hold automatically and do not need further verification.

  Choose a homogeneous basis $\mathcal Z=\{z_1, z_2, z_3\}$ for $L(2)$, where $|z_1| = |y_1y_2| = 2d_1$ and $|z_2| = |z_3| = d_1 + d_3$. In this basis
  \[ X =
    \begin{bmatrix}
      c_{\{1,2\}}  & 0 \\ 0  & X_0
    \end{bmatrix},
  \]
  for some $2 \times 2$ integer matrix $X_0$ of nonzero determinant.
  
  A homogeneous change of variables in $L(1)$ is given by a block matrix $
  \begin{bmatrix}
    A  & 0 \\ 0 & \pm 1
  \end{bmatrix}$, where $A \in \GL(2; \Z)$. The effect of this change of variables in $P$ is to change the basis by the block matrix $
  \begin{bmatrix}
    \det(A) & 0 \\ 0 & \pm A 
  \end{bmatrix}$. Therefore, by change of basis in $L(1)$ and $L(2)$ we can ensure that $X$ becomes
  \[
    X' = \begin{bmatrix}
      \pm 1 & 0 \\0  & B
    \end{bmatrix}
    \begin{bmatrix}
      c_{\{1,2\}} & 0 \\ 0 & X_0
    \end{bmatrix}
    \begin{bmatrix}
      \det(A) & 0 \\ 0 & \pm A
    \end{bmatrix}
    \]
  where $B \in \GL(2;\Z)$. We can choose $A$ and $B$ in order to bring $X_0$ to Smith normal form $
  \begin{bmatrix}
    c_{\{2,3\}} & 0 \\ 0 & c_{\{3,1\}}
  \end{bmatrix}$, and we can choose signs so that
  \[ X' =
    \begin{bmatrix}
     c_{\{1,2\}} & 0 & 0 \\ 0 & c_{\{2,3\}} & 0 \\ 0 & 0 & c_{\{3,1\}}
    \end{bmatrix} \quad \text{ where $c_{\{1,2\}}, c_{\{2,3\}}, c_{\{3,1\}} > 0$,} \]
  as required.

\item[The $d_i$ are all distinct] In this case, $L(1)$ and $L(2)$ decompose uniquely as direct sums of rank-$1$ graded summands and there are, up to choices of sign and ordering, unique homogeneous bases for each. The matrix of $X$ takes the form
  \[
    \begin{bmatrix}
      c_{\{1,2\}} & 0 & 0 \\ 0 & c_{\{2,3\}} & 0 \\ 0 & 0 & c_{\{3,1\}}
    \end{bmatrix},
  \]
  and after some choices of sign we can ensure that the $c_{\{i,j\}}$ are positive.
\end{description}

Finally, let us consider $L(3)$, which is an abelian group of rank $1$. For degree reasons, $b_1b_2b_3$ must lie in $L(3)$, and so there exists some positive integer $c_{\{1,2,3\}}$ for which $L(3)$ is generated by $\frac{1}{c_{\{1,2,3\}}} b_1b_2b_3$.
\end{proof} 

We close this section with an example illustrating why Theorem~\ref{th:ordersAllBasisG} excludes the case of two generators of the same even degree.

\begin{example} \label{ex:bad3}
  It is possible for an order in a sphere-product $\Q$-algebra not to be a weighted sphere-product ring. Observe that a weighted sphere-product ring has a minimal set of homogeneous generators $\{a_{1},\ldots,a_{m}\}$ that satisfy $a_i^2 = 0$, and all homogeneous components have bases that consist of fractional multiples of monomials in these generators. 
 In the algebra $F_\Q[x_2,y_2,z_3]/(x_2^2,y_2^2,z_3^2)$, in which $|x_2| = |y_2| = 2$ and $|z_3| = 3$, we can take the lattice $A$ generated by
  \[ 1;\, x_2, y_2;\, z_3;\, x_2y_2; \, x_2z_3, \frac{1}{2}(x_2z_3+y_2z_3); \, \frac{1}{2}x_2y_2z_3. \]
  This is closed under multiplication, and therefore constitutes an order. For $t\geq 1$, let $A_{t}$ be the submodule of $A$ consisting of elements of degree $t$. Since $x_2, y_2$ are of even degree, any basis for~$A_2$ that consists of elements that square to $0$ must be of the form $\{\pm x_2, \pm y_2\}$. Up to sign, only possible basis for $A_3$ is $\{z_3\}$. On the other hand, there is no basis of $A_5$ that consists of fractional multiples of the monomials $x_2z_3$ and $y_2z_3$.
\end{example}
Example \ref{ex:bad3} arises because the element $x_2 + y_2$, which one might wish to use as part of the distinguished generating set for a weighted sphere-product algebra is inadmissible for this purpose, since $(x_2+y_2)^2=2x_2y_2\neq 0$.

%\williams{Weighted sphere-product rings are defined in Definition \ref{def:mgwspr}. There has to be a particular set of elements that act like generators up to the introduction of factors $\frac{1}{c_\sigma}$ and, crucially, \emph{square to $0$}.

%  The order whose $\Z$-basis is
%  \[ 1;\, x, y;\, z;\, xy; \, xz, \frac{1}{2}(xz+yz); \, \frac{1}{2}xyz. \]
%  (I am being agnostic about the gradings) can be written in a different way. Let us define $u:=x+y$. Then the $\Z$-basis is
%  \[ 1; \, x, u;\, z;\, xu;\, xz, \frac{1}{2}uz; \, \frac{1}{2}xuz. \]
%  To see this, observe that $xu=x^2+xy = xy$ since $x^2=0$.
%  If $x,y$ have odd degree, then $u^2=(x+y)^2 = x^2 +xy+yx + y^2= 0 + xy-xy + 0=0$ by anticommutativity of odd-degree multiplication, and there is no problem, this is a weighted sphere-product ring using the generators $x,u,z$. If $x,y$ have even degree, $u^2=(x+y)^2 = x^2 + 2xy + y^2 = 2xy \neq 0$, so this is not a presentation as such a ring. It takes a tiny bit more effort to see that no such presentation is available.}

\appendix
\section{Alternating representations}\label{sec:altern-repr}

If $n$ is a positive integer and $V$ is a module over a commutative ring $S$, then we can form the $n$-fold \emph{alternating power} $\Alt^n_S V$, which is the quotient of the $n$-fold tensor power $V^{\tensor n}$ by the ideal generated by elementary tensors $v_1 \tensor \cdots \tensor v_n$ having repeated elements. 

If $V$ is finitely generated and free, with an ordered basis $\{e_1, \dots, e_r\}$, then $\Alt^n_S V$ is also finitely generated and free, having as a basis the images of elementary tensors $e_{i_1} \tensor \cdots \tensor e_{i_n}$ satisfying $i_1 < i_2 < \dots < i_n$. We will order this new basis lexicographically.

The construction $V \mapsto \Alt^n V$ is functorial. Note that we omit the $S$ from the notation when no confusion may arise. We obtain a group homomorphism
\[ \Alt^n : \GL(V) \to \GL(\Alt^n V). \]
When $V=\Z^r$, by virtue of our ordering convention, we may write this as
\[ \Alt^n : \GL(r; \Z) \to \GL\left(\binom{r}{n}; \Z\right). \]

\begin{proposition}\label{pr:Zn=2r=3homo}
  The homomorphism
  \[ \Alt^2 : \GL(3;\Z) \to \GL(3; \Z) \]
  has image $\SL(3;\Z)$.
\end{proposition}
\begin{proof}
  First of all, the function $g \mapsto \det(\Alt^2(g))$ is the restriction to $\GL(3;\Z)$ of a homomorphism $\GL(3;\Q) \to \Q^\times$. It is therefore given by $\det(g)^r$ for some integer $r$, and by considering a scalar matrix $g$, for instance, we deduce that
  \[ \det(\Alt^2(g)) = \det(g)^2. \]
  In particular, the image of $\Alt^2 : \GL(3; \Z) \to \GL(3; \Z)$ lies in $\SL(3; \Z)$, since every matrix in $\GL(3;\Z)$ has determinant $1$ or $-1$.

  To see that every matrix in $\SL(3; \Z)$ actually arises, we consider $\Alt^2$ applied to an elementary matrix. For example, a calculation shows
  \[\Alt^2\left(
      \begin{bmatrix}
        1 & -a & 0 \\ 0& 1 & 0 \\ 0 & 0 & 1
      \end{bmatrix} \right) =
    \begin{bmatrix}
      1 & 0 & a \\ 0 & 1 & 0 \\ 0 & 0 & 1
    \end{bmatrix}.
  \]
  The case of other elementary matrices is similar.
  
  We deduce that all elementary matrices are in the image of $\Alt^2$. A result of Vaserstein, \mbox{\cite[Corollary 3.3]{V}}, implies that elementary matrices generate $\SL(3;\Z)$.
\end{proof}

\bibliographystyle{amsalpha }

\end{document}